# General Integral Inequalities Including Weight Functions

Qian Feng and Sing Kiong Nguang

### Abstract

In this note, we present two general classes of integral inequalities motivated by their applications to infinite-dimensional systems. The inequalities possess general structures in terms of their weight functions, and the integral kernels of lower quadratic bounds which belong to a weighted $\mathbb{L}^2$ space. Many existing inequalities in the published literature, including those with free matrix variables, are the special cases of the proposed inequalities. The optimality and a relation concerning the lower bounds of the proposed inequalities are also investigated rigorously. For specific applications, our inequalities are applied to construct a Liapunov-Krasovskii functional for the stability analysis of a linear coupled differential-difference system with a distributed delay, which gives to equivalent stability conditions based on the properties of the proposed inequalities. Apart from the application demonstrated in this note, the proposed inequalities can be applied in general contexts such as the stability analysis of PDE-related systems or sampled-data systems.

## I. Introduction

Many control and optimization problems involve the applications of integral inequalities. Notable examples can be found in the stability analysis and stabilization of linear delay systems [1], [2] and PDE-related systems [3]–[6] when the direct Liapunov method is utilized. Since utilizing what kind of integral inequalities for the construction of a Liapunov functional can directly affect the conservatism of the resulting stability or stabilization conditions, thus it is of cardinal importance to develop integral inequalities with optimal mathematical structure.

The aim of this note is to present two general classes of integral inequalities following the research trajectories in [1], [2], [7]–[12] and [13]–[15], where integral inequalities were gradually generalized to construct Liapunov Krasovskii functionals with increasingly general structures. In Section II, we present the first class of integral inequalities which has no extra matrix variables other than the original matrix term in a quadratic form. Unlike the existing results in [1], [2], [7]–[12] where the integral kernels or the weight functions belong to a restrictive class of functions, the integral kernels of our inequalities can be any function belong to a weighted $\mathbb{L}^2$ space with a general weight function. This also means that many existing inequalities in [1], [2], [7]–[12] are the particular cases of the proposed inequalities. Moreover, it is proved in this note that the optimality of the first class of inequalities can be guaranteed by the principle of least-squares approximations in Hilbert space, though the integral kernels do not need to be mutually orthogonal.

The authors are with School of Electrical and Computer Engineering The University of Auckland Auckland, New Zealand. Emails: qfen204@aucklanduni.ac.nz, nguang@auckland.ac.nz

On the other hand, the second class of integral inequalities, called Free Matrix Type, is derived in Section III, which can generalize the existing inequalities in [13]–[15]. This type of inequalities can be useful in dealing with time-vary delays or sampled-data systems as mentioned in Remark 7 of [14] and Section IV. B of [16], respectively. We then prove an interesting result concerning the equivalence between the lower bounds of the two classes of inequalities under certain conditions. By this equivalence result, one can conclude that the lower bounds of many existing inequalities are essentially equivalent. To show a specific application of our inequalities, we apply them in Section IV to derive stability conditions for a linear coupled differential-difference system (CDDS) [17] with a distributed delay by constructing a parameterized complete Liapunov-Krasovskii functional. We show that equivalent stability conditions can be obtained by the application of our inequalities. The core contributions in this note are rooted in the generality of the proposed inequalities with the guarantee of optimality. This provides great potential to apply them to tackle problems in the contexts of control and optimizations.

## NOTATION

Let $\mathbb{S}^n := \{X \in \mathbb{R}^{n\times n} : X = X^\top\}$ and $\mathbb{R}_{\geq a} := \{x \in \mathbb{R} : x \geq a\}$ and $\mathbb{R}^{n\times m}_{[r]} = \{X \in \mathbb{R}^{n\times m} : \mathrm{rank}(X) = r\}$. $\mathbb{L}_f(\mathcal{X} ; \mathcal{Y})$ denotes the space of all functions which are Lebesgue integrable from the set $\mathcal{X}$ onto $\mathcal{Y}$. We frequently utilize the notations of universal quantifier $\forall$ and the existential quantifier $\exists$ in this paper. $\mathbf{Sy}(X) := X + X^\top$ stands for the sum of a matrix with its transpose. $\mathbf{Col}_{i=1}^n x_i := \left[\mathbf{Row}_{i=1}^n x_i^\top\right]^\top = \left[x_1^\top \cdots x_i^\top \cdots x_n^\top\right]^\top$ stands for a column vector containing a sequence of mathematical objects (scalars, vectors, matrices etc.). The symbol $*$ is used to indicate $[*]YX = X^\top YX$ or $X^\top Y[*] = X^\top YX$ or $\left[\begin{smallmatrix} A & B \\ * & C \end{smallmatrix}\right] = \left[\begin{smallmatrix} A & B \\ B^\top & C \end{smallmatrix}\right]$. Note that $*$ here can only represent one matrix symbol, so $[*]YXZ = Z^\top X \top YXZ$ is incorrect and $[*]Y[XZ] = Z^\top X \top YXZ$ is correct since the bracket of $[XY]$ makes it as one symbol for the matrices $XY$. Moreover, $\mathsf{O}_{n\times n}$ denotes a $n \times n$ zero matrix which can be abbreviated into $\mathsf{O}_n$, while $\mathbf{0}_n$ represents a $n \times 1$ column vector. We frequently use $X \oplus Y = \left[\begin{smallmatrix} X & \mathsf{O} \\ * & Y \end{smallmatrix}\right]$ to denote the diagonal sum of two matrices, respectively. $\otimes$ stands for the Kronecker product. Finally, we assume the order of matrix operations to be *matrix (scalars) multiplications* $> \otimes > \oplus > +$.

## II. FIRST CLASS OF INTEGRAL INEQUALITIES

In this section, we present the first general class of inequalities without introducing slack matrix variables. We will frequently apply the following lemma throughout the entire paper, which is a particular case of the property of the Kronecker product $(A \otimes B)(C \otimes D)$.

**Lemma 1.** $\forall X \in \mathbb{R}^{n\times m},\ \forall Y \in \mathbb{R}^{m\times p},\ \ \forall Z \in \mathbb{R}^{q\times r}$,

$$(X \otimes I_q)(Y \otimes Z) = (XY) \otimes (I_q Z) = (XY) \otimes Z = (XY) \otimes (ZI_r) = (X \otimes Z)(Y \otimes I_r). \tag{1}$$

The following lemma is partially taken from Lemma 4.1 in [16], which is vital for the derivation of the results in this note.

**Lemma 2.** *Let $C \in \mathbb{S}^m_{\succ 0}$ and $B \in \mathbb{R}^{m\times n}$ be given, then*

$$\forall M \in \mathbb{R}^{m\times n},\ \ B^\top C^{-1} B \succeq M^\top B + B^\top M - M^\top CM \tag{2}$$

*where the inequality in* (2) *becomes an equality with $M = C^{-1}B$.*

To present the first class of inequalities in this section, we define the following weighted Lebesgue function space

$$\mathbb{L}^2_\varpi(\mathcal{K} ; \mathbb{R}^d) := \left\{\phi(\cdot) \in \mathbb{L}_f(\mathcal{K} ; \mathbb{R}^d) : \|\phi(\cdot)\|_{2,\varpi} < \infty\right\} \tag{3}$$

with $d \in \mathbb{N}$ and $\|\phi(\cdot)\|_{2,\varpi}^2 := \int_{\mathcal{K}} \varpi(\tau)\phi^\top(\tau)\phi(\tau)\mathsf{d}\tau$, where $\varpi(\cdot) \in \mathbb{L}_f(\mathcal{K} ⨾ \mathbb{R}_{\geq 0})$ and $\varpi(\cdot)$ has only countably infinite or finite numbers of zero values. Furthermore, $\mathcal{K}$ in (3) satisfies $\mathcal{K} \subseteq \mathbb{R} \cup \{\pm\infty\}$ and the Lebesgue measure of $\mathcal{K}$ is non-zero.

**Theorem 1.** *Let $\varpi(\cdot)$, $\mathcal{K}$ and $d \in \mathbb{N}$ in* (3) *be given, and assume $\boldsymbol{f}(\cdot) = \mathbf{Col}_{i=1}^d f_i(\cdot) \in \mathbb{L}_\varpi^2\left(\mathcal{K} ⨾ \mathbb{R}^d\right)$ which satisfies*

$$\int_{\mathcal{K}} \varpi(\tau)\boldsymbol{f}(\tau)\boldsymbol{f}^\top(\tau)\mathsf{d}\tau \succ 0, \tag{4}$$

*then the inequality*

$$\int_{\mathcal{K}} \varpi(\tau)\boldsymbol{x}^\top(\tau)U\boldsymbol{x}(\tau)\mathsf{d}\tau \geq \int_{\mathcal{K}} \varpi(\tau)\boldsymbol{x}^\top(\tau)F^\top(\tau)\mathsf{d}\tau\,(\mathsf{F}\otimes U)\int_{\mathcal{K}} \varpi(\tau)F(\tau)\boldsymbol{x}(\tau)\mathsf{d}\tau \tag{5}$$

*holds for all $\boldsymbol{x}(\cdot) \in \mathbb{L}_\varpi^2(\mathcal{K} ⨾ \mathbb{R}^n)$ and $U \in \mathbb{S}_{\succeq 0}^n$, where $F(\tau) = \boldsymbol{f}(\tau) \otimes I_n$ and $\mathsf{F}^{-1} = \int_{\mathcal{K}} \varpi(\tau)\boldsymbol{f}(\tau)\boldsymbol{f}^\top(\tau)\mathsf{d}\tau$. Moreover, the optimality of* (5) *with $U \succ 0$ is guaranteed by the application of* (2)*.*

*Proof:* Let $\boldsymbol{\varepsilon}(\tau) = \boldsymbol{x}(\tau) - F^\top(\tau)\boldsymbol{\omega}$ with $\boldsymbol{\omega} \in \mathbb{R}^{dn}$ and $F(\tau) = \boldsymbol{f}(\tau) \otimes I_n$, through which we have

$$\int_{\mathcal{K}} \varpi(\tau)\boldsymbol{\varepsilon}^\top(\tau)U\boldsymbol{\varepsilon}(\tau)\mathsf{d}\tau = \int_{\mathcal{K}} \varpi(\tau)\boldsymbol{x}^\top(\tau)U\boldsymbol{x}(\tau)\mathsf{d}\tau - 2\int_{\mathcal{K}} \varpi(\tau)\boldsymbol{x}^\top(\tau)UF^\top(\tau)\mathsf{d}\tau\boldsymbol{\omega} + \boldsymbol{\omega}^\top\int_{\mathcal{K}} \varpi(\tau)F(\tau)UF^\top(\tau)\mathsf{d}\tau\boldsymbol{\omega}. \tag{6}$$

Now apply (1) to $UF^\top(\tau)$ in (6) with $F(\tau) = \boldsymbol{f}(\tau) \otimes I_n$, we have

$$UF^\top(\tau) = U(\boldsymbol{f}^\top(\tau) \otimes I_n) = \boldsymbol{f}^\top(\tau) \otimes U = F^\top(\tau)(I_d \otimes U). \tag{7}$$

Moreover, applying (7) to $\int_{\mathcal{K}} \varpi(\tau)\boldsymbol{x}^\top(\tau)UF^\top(\tau)\mathsf{d}\tau\boldsymbol{\omega}$ in (6) yields

$$\int_{\mathcal{K}} \varpi(\tau)\boldsymbol{x}^\top(\tau)UF^\top(\tau)\mathsf{d}\tau\boldsymbol{\omega} = \int_{\mathcal{K}} \varpi(\tau)\boldsymbol{x}^\top(\tau)F^\top(\tau)\mathsf{d}\tau\,(I_d\otimes U)\,\boldsymbol{\omega} = \boldsymbol{\vartheta}^\top(I_d\otimes U)\boldsymbol{\omega} \tag{8}$$

with $\boldsymbol{\vartheta} = \int_{\mathcal{K}} \varpi(\tau)F(\tau)\boldsymbol{x}(\tau)\mathsf{d}\tau$. By (7) and (1), we have

$$\int_{\mathcal{K}} \varpi(\tau)F(\tau)UF^\top(\tau)\mathsf{d}\tau = \int_{\mathcal{K}} \varpi(\tau)F(\tau)F^\top(\tau)\mathsf{d}\tau(I_d\otimes U) = \left[\int_{\mathcal{K}} \varpi(\tau)\boldsymbol{f}(\tau)\boldsymbol{f}^\top(\tau)\mathsf{d}\tau \otimes I_n\right](I_d\otimes U) = \mathsf{F}^{-1}\otimes U \tag{9}$$

where $\mathsf{F}^{-1} = \int_{\mathcal{K}} \varpi(\tau)\boldsymbol{f}(\tau)\boldsymbol{f}^\top(\tau)\mathsf{d}\tau \succ 0$ given (4). Moreover, the expressions in (6) can be simplified into

$$\int_{\mathcal{K}} \varpi(\tau)\boldsymbol{\varepsilon}^\top(\tau)U\boldsymbol{\varepsilon}(\tau)\mathsf{d}\tau = \int_{\mathcal{K}} \varpi(\tau)\boldsymbol{x}^\top(\tau)U\boldsymbol{x}(\tau)\mathsf{d}\tau - \mathbf{Sy}\left(\boldsymbol{\vartheta}^\top\left(I_d\otimes U\right)\boldsymbol{\omega}\right) + \boldsymbol{\omega}^\top\left(\mathsf{F}^{-1}\otimes U\right)\boldsymbol{\omega}. \tag{10}$$

by (8) and (9). Since $\int_{\mathcal{K}} \varpi(\tau)\boldsymbol{\varepsilon}^\top(\tau)U\boldsymbol{\varepsilon}(\tau)\mathsf{d}\tau \geq 0$ if $U \succeq 0$ in (10), it yields

$$\int_{\mathcal{K}} \varpi(\tau)\boldsymbol{x}^\top(\tau)U\boldsymbol{x}(\tau)\mathsf{d}\tau \geq \mathbf{Sy}\left(\boldsymbol{\vartheta}^\top\left(\mathsf{F}^{-1}\otimes U\right)\boldsymbol{\omega}\right) - \boldsymbol{\omega}^\top\left(\mathsf{F}\otimes U\right)\boldsymbol{\omega} \tag{11}$$

for all $\boldsymbol{x}(\cdot) \in \mathbb{L}_\varpi^2(\mathcal{K} ⨾ \mathbb{R}^n)$ and $U \in \mathbb{S}_{\succeq 0}^n$. Moreover, let $\boldsymbol{\omega} = (\mathsf{F}\otimes I_n)\boldsymbol{\vartheta}$, then (10) gives (5), thus (5) is proved.

Now we proceed to show the optimality of (5) with $U \succ 0$ via the application of Lemma 2. Applying (2) to (11) with $C = \mathsf{F}^{-1}\otimes U$, $B = (I_d\otimes U)\boldsymbol{\vartheta}$ and $M = \boldsymbol{\omega}$, one can conclude that

$$\mathbf{Sy}\left(\boldsymbol{\vartheta}^\top\left(I_n\otimes U\right)\boldsymbol{\omega}\right) - \boldsymbol{\omega}^\top\left(\mathsf{F}^{-1}\otimes U\right)\boldsymbol{\omega} \leq \boldsymbol{\vartheta}^\top\left(I_d\otimes U\right)\left(\mathsf{F}\otimes U^{-1}\right)\left(I_d\otimes U\right)\boldsymbol{\vartheta} = \boldsymbol{\vartheta}^\top(\mathsf{F}\otimes U)\boldsymbol{\vartheta} \tag{12}$$

for all $\boldsymbol{x}(\cdot) \in \mathbb{L}_\varpi^2(\mathcal{K} ⨾ \mathbb{R}^n)$ and $\boldsymbol{\omega} \in \mathbb{R}^{dn}$, where the inequality in (12) becomes an equality with $M = \boldsymbol{\omega} = C^{-1}B = (\mathsf{F}\otimes I_n)\boldsymbol{\vartheta}$. Thus all the statements in Theorem 1 are proved. ■

*Remark* 1. By considering Theorem 7.2.10 in [18], the condition in (4) indicates that the functions $\{f_i(\cdot)\}_{i=1}^d$ in $\boldsymbol{f}(\cdot) = \mathbf{Col}_{i=1}^d f_i(\tau)$ are linearly independent in a Lebesgue sense. Given the general mathematical structure of $\boldsymbol{f}(\cdot) \in \mathbb{L}^2_{\varpi}(\mathcal{K}⨟\mathbb{R}^d)$ and $\varpi(\cdot)$ in Theorem 1, the generality of (5) is evident. For instance, $\boldsymbol{f}(\cdot)$ can contain orthogonal functions [8], elementary functions or other types of functions as long as $\{f_i(\cdot)\}_{i=1}^d$ in $\boldsymbol{f}(\cdot) = \mathbf{Col}_{i=1}^d f_i(\tau)$ are linearly independent in a Lebesgue sense.

The optimality of (5) is inferred in Theorem 1 via the application of the matrix inequality in Lemma 2. In the following corollary, we show that the optimality of (5) with $U \succ 0$ is also ensured by the least-squares approximations in Hilbert space, which coincides with the result produced via Lemma 2.

**Corollary 1.** *With $U \succ 0$, the optimality of (5) is guaranteed by the least-squares approximations in $\mathbb{L}^2_{\varpi}(\mathcal{K}⨟\mathbb{R})$, which is consistent with the result produced via Lemma 2 in the proof of Theorem 1.*

*Proof.* By the definitions of matrix multiplication and the Kronecker product, $\boldsymbol{\varepsilon}(\tau) = \boldsymbol{x}(\tau) - F^\top(\tau)\boldsymbol{\omega}$ at the beginning of the proof of Theorem 1 can be rewritten into

$$\boldsymbol{\varepsilon}(\tau) = \boldsymbol{x}(\tau) - (\boldsymbol{f}^\top(\tau) \otimes I_n)\boldsymbol{\omega} = \boldsymbol{x}(\tau) - \sum_{i=1}^d f_i(\tau)\boldsymbol{\omega}_i = \boldsymbol{x}(\tau) - \left[\mathbf{Row}_{i=1}^d \boldsymbol{\omega}_i\right] \boldsymbol{f}(\tau) = \boldsymbol{x}(\tau) - \left[\mathbf{Col}_{i=1}^n \boldsymbol{\alpha}_i^\top\right] \boldsymbol{f}(\tau)$$

where $\boldsymbol{\omega} = \mathbf{Col}_{i=1}^d \boldsymbol{\omega}_i$ with $\boldsymbol{\omega}_i \in \mathbb{R}^n$, and $\mathbf{Row}_{i=1}^d \boldsymbol{\omega}_i = \mathbf{Col}_{i=1}^n \boldsymbol{\alpha}_j^\top$ with $\boldsymbol{\alpha}_j \in \mathbb{R}^d$. Thus, it is clear that $\boldsymbol{\omega}$ can be interpreted as an approximation coefficient, and $\boldsymbol{\varepsilon}(\cdot)$ measures the error of approximation. From the structure of $\boldsymbol{\varepsilon}(\tau) = \boldsymbol{x}(t) - \left[\mathbf{Row}_{i=1}^d \boldsymbol{\omega}_i\right] \boldsymbol{f}(\tau)$, it shows that each row of $\boldsymbol{x}(\cdot) = \mathbf{Col}_{i=1}^n x_i(\tau)$ is individually approximated by all the functions $\{f_i(\cdot)\}_{i=1}^d$ in $\boldsymbol{f}(\cdot) = \mathbf{Col}_{i=1}^d f_i(\cdot)$. Consider the norm square (Inner product) $\langle \boldsymbol{\varepsilon}(\cdot), \boldsymbol{\varepsilon}(\cdot)\rangle_{\varpi} = \|\boldsymbol{\varepsilon}(\cdot)\|^2_{\varpi} := \int_{\mathcal{K}} \varpi(\tau)\boldsymbol{\varepsilon}^\top(\tau)\boldsymbol{\varepsilon}(\tau)\mathrm{d}\tau = \sum_{i=1}^n \int_{\mathcal{K}} \varpi(\tau)\varepsilon_i^2(\tau)\mathrm{d}\tau$, the objective now is to calculate $\min_{\boldsymbol{\omega}} \|\boldsymbol{\varepsilon}(\cdot)\|^2_{\varpi}$ for a given $\boldsymbol{f}(\cdot) \in \mathbb{L}^2_{\varpi}(\mathcal{K}⨟\mathbb{R}^d)$ by finding an appropriate $\boldsymbol{\omega}$.

The minimization problem can be solved by using the orthogonality principle, namely, least-squares approximations [19] in the Hilbert Space $\left(\mathbb{L}^2_{\varpi}(\mathcal{K}⨟\mathbb{R}), \langle\cdot,\cdot\rangle_{\varpi}\right)$. Namely, for any $x_i(\cdot) \in \mathbb{L}^2_{\varpi}(\mathcal{K}⨟\mathbb{R})$ we want to find the closest point $\widehat{f}_i(\cdot) = \boldsymbol{\alpha}_i^\top \boldsymbol{f}(\tau)$ to $x_i(\cdot)$ with an appropriate coefficient $\boldsymbol{\omega}$, where $\widehat{f}_i(\cdot)$ belongs to the liner subspace of $\left(\mathbb{L}^2_{\varpi}(\mathcal{K}⨟\mathbb{R}), \langle\cdot,\cdot\rangle_{\varpi}\right)$ spanned by the vectors $\{f_i(\cdot)\}_i^d \subset \mathbb{L}^2_{\varpi}(\mathcal{K}⨟\mathbb{R})$. For a rigorous presentation of the theory on this topic, see section 10.2 in [19].

According to the procedures of least-squares approximations outlined in page 182 of [19], the coefficient corresponds to the minimization of $\|\boldsymbol{\varepsilon}(\cdot)\|^2_{\varpi}$ must satisfy

$$\mathsf{F}^{-1}\boldsymbol{\alpha}_i = \int_{\mathcal{K}} \varpi(\tau)x_i(\tau)\boldsymbol{f}(\tau)\mathrm{d}\tau \in \mathbb{R}^d, \quad i = 1 \cdots n$$

which can be further formulated into the compact form

$$\mathsf{F}^{-1} \underset{i=1}{\overset{n}{\mathbf{Row}}}\, \boldsymbol{\alpha}_i = \mathsf{F}^{-1} \underset{i=1}{\overset{d}{\mathbf{Col}}}\, \boldsymbol{\omega}_i^\top = \int_{\mathcal{K}} \varpi(\tau)\boldsymbol{f}(\tau)\boldsymbol{x}^\top(\tau)\mathrm{d}\tau \in \mathbb{R}^{d\times n}.$$

Thus we have $\mathbf{Row}_{i=1}^d \boldsymbol{\omega}_i = \int_{\mathcal{K}} \varpi(\tau)\boldsymbol{x}(\tau)\boldsymbol{f}^\top(\tau)\mathrm{d}\tau\mathsf{F}$. Now by the definition of matrix multiplication with $\widetilde{\boldsymbol{f}}(\tau) := \mathbf{Col}_{i=1}^d \widetilde{f}(\tau) := \mathsf{F}\boldsymbol{f}(\tau)$, the previous expression can be changed into $\mathbf{Row}_{i=1}^d \boldsymbol{\omega}_i = \int_{\mathcal{K}} \varpi(\tau) \mathbf{Row}_{i=1}^d \boldsymbol{x}(\tau)\widetilde{f}_i(\tau)\mathrm{d}\tau$. Hence we can conclude that

$$\boldsymbol{\omega}_i = \int_{\mathcal{K}} \varpi(\tau)\boldsymbol{x}(\tau)\widetilde{f}_i(\tau)\mathrm{d}\tau, \quad \forall i = 1 \cdots d.$$

Given $\boldsymbol{\omega} = \mathbf{Col}_{i=1}^d \boldsymbol{\omega}_i \in \mathbb{R}^{dn}$ as we had defined at the beginning of the proof, the above equation can be further formulated into

$$\boldsymbol{\omega} = \underset{i=1}{\overset{d}{\mathsf{Col}}} \int_{\mathcal{K}} \varpi(\tau)\widetilde{f}_i(\tau)\boldsymbol{x}(\tau)\mathsf{d}\tau = \int_{\mathcal{K}} \varpi(\tau)(\widetilde{\boldsymbol{f}}(\tau)\otimes I_n)\boldsymbol{x}(\tau)\mathsf{d}\tau$$

$$= (\mathsf{F}\otimes I_n)\int_{\mathcal{K}} \varpi(\tau)(\boldsymbol{f}(\tau)\otimes I_n)\boldsymbol{x}(\tau)\mathsf{d}\tau = (\mathsf{F}\otimes I_n)\boldsymbol{\vartheta} \quad (13)$$

considering the fact that $\widetilde{\boldsymbol{f}}(\tau) = \mathsf{F}\boldsymbol{f}(\tau)$. In the proof of Theorem 1, we have shown that the largest inequality's lower bound of (11) is given by $\boldsymbol{\omega} = (\mathsf{F}\otimes I_n)\boldsymbol{\vartheta}$ leading to (5). This indicates that the choice of $\boldsymbol{\omega} = (\mathsf{F}\otimes I_n)\boldsymbol{\vartheta}$ made for (11) coincides with the coefficient in (13) calculated based on the principle of least-squares approximations which minimizes $\|\boldsymbol{\varepsilon}(\cdot)\|^2$. Furthermore, since for all $U \succ 0$ there exists $N \in \mathbb{R}^{n\times n}_{[n]}$ such that $U = N^\top N$, then the structure of $\int_{\mathcal{K}} \varpi(\tau)\boldsymbol{\varepsilon}^\top(\tau)U\boldsymbol{\varepsilon}(\tau)\mathsf{d}\tau = \int_{\mathcal{K}} \varpi(\tau)\boldsymbol{\varepsilon}^\top(\tau)N^\top N\boldsymbol{\varepsilon}(\tau)\mathsf{d}\tau = \int_{\mathcal{K}} \varpi(\tau)\widetilde{\boldsymbol{\varepsilon}}^\top(\tau)\widetilde{\boldsymbol{\varepsilon}}(\tau)\mathsf{d}\tau = \sum_{i=1}^n \int_{\mathcal{K}} \varpi(\tau)\widetilde{\varepsilon}_i^2(\tau)\mathsf{d}\tau$ is equivalent to $\int_{\mathcal{K}} \varpi(\tau)\boldsymbol{\varepsilon}^\top(\tau)\boldsymbol{\varepsilon}(\tau)\mathsf{d}\tau = \int_{\mathcal{K}} \varpi(\tau)\varpi(\tau)\varepsilon_i^2(\tau)\mathsf{d}\tau$ for $\mathbb{L}^2_\varpi(\mathcal{K} \,; \mathbb{R})$. Consequently, the coefficient in (13) calculated via least-squares approximations also minimizes $\int_{\mathcal{K}} \varpi(\tau)\boldsymbol{\varepsilon}^\top(\tau)U\boldsymbol{\varepsilon}(\tau)\mathsf{d}\tau$ in (6). ■

*Remark* 2. The optimal value of $\boldsymbol{\omega}$ for the minimization of $\int_{\mathcal{K}} \varpi(\tau)\boldsymbol{\varepsilon}^\top(\tau)U\boldsymbol{\varepsilon}(\tau)\mathsf{d}\tau$ may also be determined by differentiating $\boldsymbol{\omega}$ in (10). This kind of idea has been considered in page 2 of [20] and the proof of Lemma 3 in [21], which may also require the application of Hessian [21]. On the other hand, the derivation of the optimal value of $\boldsymbol{\omega}$ is very straightforward via the application of Lemma 2 at (12), which is consistent with the result produced by the principle of least-squares approximations in Corollary 1.

*Remark* 3. With appropriate $\boldsymbol{x}(\cdot)$, $\mathcal{K} = [-r, 0]$ and mathematical manipulations, the inequalities in [1, eq. (5)–(6)] can be obtained by (5) with $\boldsymbol{f}(\cdot)$ containing the Legendre polynomials over $[-r, 0]$. Meanwhile, if $\boldsymbol{f}(\cdot)$ contains only orthogonal functions, then [8, eq.(5)] can be obtained by (5) via the use of commutation matrices [22]. Let $\mathcal{K} = [0, +\infty]$, then the inequality in [12, eq. (9)] is the special case of (5) with appropriate $\varpi(\cdot)$ and $\boldsymbol{f}(\cdot)$. Finally, with $\varpi(\tau) = 1$, (5) becomes [2, eq. (16)].

*Remark* 4. One can conclude that the polynomials in [9, eq. (13)–(14)]; [11, eq. (3)–(4)]; [10, eq.(2)] and [21, eq.(5)-(7)] are the special cases of the Jacobi polynomials [23, See 22.3.2] over $[a, b]$ with appropriate weight functions. Note that the above conclusions on [21, eq.(5)-(7)] and [11, eq. (3)–(4)] are established based on the use of the Cauchy formula [24, page 193] for repeated integrations [9, See eq.(5)–(6) and eq.(25)–(26)]. As a result, the inequalities in [9, eq.(27),(34)]; [10, eq.(2)]; [11, eq.(1)–(2)]; and the inequalities in [21, (8)–(9)] with finite terms of summations, can be obtained by (5) with appropriate $\boldsymbol{x}(\cdot)$ and choosing $\boldsymbol{f}(\cdot)$ be the Jacobi polynomials over $[a, b]$ with the corresponding $\varpi(\tau) = (\tau - a)^p$ or $\varpi(\tau) = (b - \tau)^p$. Finally, see the discussions in [25, Pages 50-53] where the inclusion of the existing inequalities by (5) is proved with more mathematical details.

*Remark* 5. Note that a substitution $G\boldsymbol{f}(\tau) \to \boldsymbol{f}(\tau)$ for (5) with an invertible $G \in \mathbb{R}^{n\times n}$ gives a lower bound which is equivalent to the inequality's lower bound of (5).

## III. Second Class of integral Inequalities

This section is devoted to present another general class of inequalities named as the free matrix type. This type of inequalities has been previously researched in [13]–[15], which can be useful in dealing with the stability analysis of systems with time-varying delays. Finally, a relation between the proposed inequalities in Sections II and III is established in a theorem.

The following lemma, which will be applied for the derivations of Theorem 2, can be obtained via the definition of matrix multiplication and the Kronecker product.

**Lemma 3.** *Given a matrix* $X := \mathbf{Row}_{i=1}^{d} X_i \in \mathbb{R}^{n\times d\rho n}$ *with* $n; d; \rho \in \mathbb{N}$ *and a function* $\boldsymbol{f}(\tau) = \mathbf{Col}_{i=1}^{d} f_i(\tau) \in \mathbb{R}^d$, *we have*

$$X(\boldsymbol{f}(\tau)\otimes I_{\rho n}) = \sum_{i=1}^{d} f_i(\tau)X_i = \left(\boldsymbol{f}^\top(\tau)\otimes I_n\right)\widehat{X} \tag{14}$$

*where* $\widehat{X} := \mathbf{Col}_{i=1}^{d} X_i \in \mathbb{R}^{dn\times\rho n}$.

**Theorem 2.** *Given* $\rho \in \mathbb{N}$ *and the same* $\varpi(\cdot)$, $\boldsymbol{f}(\cdot)$ *and* $U$ *in Theorem 1, if there exist* $Y \in \mathbb{S}^{\rho dn}$ *and* $X = \mathbf{Row}_{i=1}^{d} X_i \in \mathbb{R}^{n\times\rho dn}$ *such that*

$$\begin{bmatrix} U & -X \\ * & Y \end{bmatrix} \succeq 0, \tag{15}$$

*then the inequality*

$$\int_{\mathcal{K}} \varpi(\tau)\boldsymbol{x}^\top(\tau)U\boldsymbol{x}(\tau)\mathsf{d}\tau \geq \mathbf{Sy}\left(\boldsymbol{\vartheta}^\top\widehat{X}\boldsymbol{z}\right) - \boldsymbol{z}^\top\mathsf{W}\boldsymbol{z} \tag{16}$$

*holds for all* $\boldsymbol{z} \in \mathbb{R}^{\rho n}$ *and* $\boldsymbol{x}(\cdot) \in \mathbb{L}^2_{\varpi}(\mathcal{K} ⨾ \mathbb{R}^n)$, *where* $\boldsymbol{\vartheta} = \int_{\mathcal{K}} \varpi(\tau)F(\tau)\boldsymbol{x}(\tau)\mathsf{d}\tau$ *and* $\widehat{X} = \mathbf{Col}_{i=1}^{d} X_i \in \mathbb{R}^{dn\times\rho n}$ *and* $\mathsf{W} = \int_{\mathcal{K}} \varpi(\tau)(\boldsymbol{f}^\top(\tau)\otimes I_{\rho n})Y(\boldsymbol{f}(\tau)\otimes I_{\rho n})\mathsf{d}\tau \in \mathbb{S}^{\rho n}$. *Furthermore, let* $\Upsilon \in \mathbb{R}^{dn\times\rho n}$ *and* $\boldsymbol{z} \in \mathbb{R}^{\rho n}$ *be any matrix and vector satisfy*

$$\Upsilon\boldsymbol{z} = \int_{\mathcal{K}} \varpi(\tau)F(\tau)\boldsymbol{x}(\tau)\mathsf{d}\tau \tag{17}$$

*for all* $\boldsymbol{x}(\cdot) \in \mathbb{L}^2_{\varpi}(\mathcal{K} ⨾ \mathbb{R}^n)$, *then we have*

$$\int_{\mathcal{K}} \varpi(\tau)\boldsymbol{x}^\top(\tau)U\boldsymbol{x}(\tau)\mathsf{d}\tau \geq \boldsymbol{z}^\top\left[\mathbf{Sy}\left(\Upsilon^\top\widehat{X}\right) - \mathsf{W}\right]\boldsymbol{z} \tag{18}$$

*for all* $\boldsymbol{x}(\cdot) \in \mathbb{L}^2_{\varpi}(\mathcal{K} ⨾ \mathbb{R}^n)$.

*Proof:* Given (15), we have

$$\begin{aligned}
&\int_{\mathcal{K}} \varpi(\tau)\begin{bmatrix}\boldsymbol{x}(\tau)\\ \boldsymbol{f}(\tau)\otimes\boldsymbol{z}\end{bmatrix}^\top \begin{bmatrix} U & -X \\ * & Y \end{bmatrix}\begin{bmatrix}\boldsymbol{x}(\tau)\\ \boldsymbol{f}(\tau)\otimes\boldsymbol{z}\end{bmatrix}\mathsf{d}\tau \\
&\qquad = \int_{\mathcal{K}} \varpi(\tau)\boldsymbol{x}^\top(\tau)U\boldsymbol{x}(\tau)\mathsf{d}\tau - \mathbf{Sy}\left[\int_{\mathcal{K}} \varpi(\tau)\boldsymbol{x}^\top(\tau)X\left(\boldsymbol{f}(\tau)\otimes\boldsymbol{z}\right)\mathsf{d}\tau\right] \\
&\qquad\quad + \int_{\mathcal{K}} \varpi(\tau)\left(\boldsymbol{f}(\tau)\otimes\boldsymbol{z}\right)^\top Y\left(\boldsymbol{f}(\tau)\otimes\boldsymbol{z}\right)\mathsf{d}\tau \geq 0.
\end{aligned} \tag{19}$$

Now using (1) and (14) to the terms in (19) yields

$$\begin{aligned}
\int_{\mathcal{K}} \varpi(\tau)\boldsymbol{x}^\top(\tau)X\left(\boldsymbol{f}(\tau)\otimes\boldsymbol{z}\right)\mathsf{d}\tau &= \int_{\mathcal{K}} \varpi(\tau)\boldsymbol{x}^\top(\tau)X(\boldsymbol{f}(\tau)\otimes I_{\rho n})\mathsf{d}\tau\boldsymbol{z} \\
&= \int_{\mathcal{K}} \varpi(\tau)\boldsymbol{x}^\top(\tau)(\boldsymbol{f}^\top(\tau)\otimes I_n)\mathsf{d}\tau\widehat{X}\boldsymbol{z} = \boldsymbol{\vartheta}^\top\widehat{X}\boldsymbol{z},
\end{aligned} \tag{20}$$

$$\int_{\mathcal{K}} \varpi(\tau)\left(\boldsymbol{f}(\tau)\otimes\boldsymbol{z}\right)^\top Y\left(\boldsymbol{f}(\tau)\otimes\boldsymbol{z}\right)\mathsf{d}\tau = \boldsymbol{z}^\top\int_{\mathcal{K}} \varpi(\tau)\left(\boldsymbol{f}^\top(\tau)\otimes I_{\rho n}\right)Y\left(\boldsymbol{f}(\tau)\otimes I_{\rho n}\right)\mathsf{d}\tau\boldsymbol{z} = \boldsymbol{z}^\top\mathsf{W}\boldsymbol{z} \tag{21}$$

where $X = \mathbf{Row}_{i=1}^{d} X_i \in \mathbb{R}^{n\times\rho dn}$ and $\widehat{X} = \mathbf{Col}_{i=1}^{d} X_i$. Substituting (20)–(21) into (19) yields (16). By considering the equality in (17) with $\Upsilon \in \mathbb{R}^{dn\times\rho n}$, then (16) becomes (18). ■

Since $\boldsymbol{f}(\cdot)$ in Theorem 2 satisfy $\boldsymbol{f}(\cdot) \in \mathbb{L}^2_{\varpi}\left(\mathcal{K} ⨾ \mathbb{R}^d\right)$ and (4) as in Theorem 1, hence the proposed free matrix type inequalities in Theorem 2 possess more general structures compared to existing results in the literature.

*Remark* 6. Let $\mathcal{K} = [a, b]$, $\varpi(\cdot) = 1$ and $\boldsymbol{f}(\tau)$ consists of the Legendre polynomials over $\mathcal{K}$, then eq.(5) in [14] can be obtained by (18) with appropriate $\Upsilon$ and $\boldsymbol{z}$ and the substitution $\boldsymbol{x}(\cdot) \leftarrow \dot{\boldsymbol{x}}(\cdot)$. Furthermore, let $\mathcal{K} = [a, b]$, $\varpi(\cdot) = 1$ and $\boldsymbol{f}(\cdot) = \mathbf{Col}_{i=0}^{m} \frac{(r-a)^i}{(b-a)^i}$, then by (16) with the substitution $\dot{\boldsymbol{x}}(t) \leftarrow \boldsymbol{x}(t)$ and considering the properties in [15, eq.(4)–(5)], one can obtain an inequality which is equivalent to eq.(7) in [15].

The following theorem shows an “equivalent" relation between the inequality’s lower bounds of (5) and (18).

**Theorem 3.** *By choosing the same $\varpi(\cdot)$, $U$ and $\boldsymbol{f}(\cdot)$ for Theorems 1 and 2 with $U \succ 0$, then the inequality’s lower bound of* (5) *corresponds to the largest inequality’s lower bound of* (18).

*Proof.* Let $U \succ 0$, $\varpi(\cdot)$ and $\boldsymbol{f}(\cdot)$ in Theorem 1 be given. Using the Schur complement [26, Theorem 1.12] with $U \succ 0$ to (15) concludes that (15) holds if and only of $Y \succeq X^\top U^{-1} X$. Now consider $\mathsf{W}$ in Theorem 2 with $Y \succeq X^\top U^{-1} X$ and (14), we have

$$\mathsf{W} \succeq \int_{\mathcal{K}} \varpi(\tau)(\boldsymbol{f}^\top(\tau) \otimes I_{nd}) X^\top U^{-1} X (\boldsymbol{f}(\tau) \otimes I_{nd}) \mathrm{d}\tau = \widehat{X}^\top \int_{\mathcal{K}} \varpi(\tau)(\boldsymbol{f}(\tau) \otimes I_n) U^{-1} \left(\boldsymbol{f}^\top(\tau) \otimes I_n\right) \mathrm{d}\tau \widehat{X}$$
$$= \widehat{X}^\top \int_{\mathcal{K}} \varpi(\tau) \left(I_d \otimes U^{-1}\right) (\boldsymbol{f}(\tau)\boldsymbol{f}^\top(\tau) \otimes I_n) \mathrm{d}\tau \widehat{X} = \widehat{X}^\top \left(\mathsf{F}^{-1} \otimes U^{-1}\right) \widehat{X} \tag{22}$$

where $\mathsf{F}^{-1} = \int_{\mathcal{K}} \varpi(\tau)\boldsymbol{f}(\tau)\boldsymbol{f}^\top(\tau)\mathrm{d}\tau \in \mathbb{S}_{\succ 0}^{d}$ and $\widehat{X} = \mathbf{Col}_{i=1}^{d} X_i \in \mathbb{R}^{dn \times \rho n}$. It is clear that $\succeq$ in (22) becomes an equality with $Y = X^\top U^{-1} X$ for any $U \succ 0$ and $X \in \mathbb{R}^{n \times \rho d n}$.

By (22) and the application of (2), we have

$$\mathbf{Sy}\left(\Upsilon^\top \widehat{X}\right) - \mathsf{W} \preceq \mathbf{Sy}\left(\Upsilon^\top \widehat{X}\right) - \widehat{X}^\top \left(\mathsf{F}^{-1} \otimes U^{-1}\right) \widehat{X} \preceq \Upsilon^\top \left(\mathsf{F} \otimes U\right) \Upsilon \tag{23}$$

holds for any $Y$ and $X$ satisfying (15) with $\Upsilon \in \mathbb{R}^{dn \times \rho n}$ in (17). Moreover, by Lemma 1, one can conclude that the two inequalities in (23) become equalities with $\widehat{X} = (\mathsf{F} \otimes U)\Upsilon$ and $Y = X^\top U^{-1} X$ where the value of $X$ here can be uniquely determined by $\widehat{X} = (\mathsf{F} \otimes U)\Upsilon$ with given $U$ and $\Upsilon$. As a result, we have

$$\boldsymbol{z}^\top \left[\mathbf{Sy}\left(\Upsilon^\top \widehat{X}\right) - \mathsf{W}\right] \boldsymbol{z} \leq \boldsymbol{z}^\top \Upsilon^\top (\mathsf{F} \otimes U)\Upsilon \boldsymbol{z} = \int_{\mathcal{K}} \varpi(\tau)\boldsymbol{x}^\top(\tau)F^\top(\tau)\mathrm{d}\tau\,(\mathsf{F} \otimes U) \int_{\mathcal{K}} \varpi(\tau)F(\tau)\boldsymbol{x}(\tau)\mathrm{d}\tau \tag{24}$$

holds for any $Y$ and $X$ satisfying (15) with $U \succ 0$, and the inequality in (24) becomes an equality with $\widehat{X} = (\mathsf{F} \otimes U)\Upsilon$ and $Y = X^\top U^{-1} X$. Hence the above arguments show that under the same $\varpi(\cdot)$, $U \succ 0$ and $\boldsymbol{f}(\cdot)$, one can always find $X$ and $Y$ for (15) to render (18) to become identical to (5) which corresponds to the largest inequality’s lower bound of (18). ■

*Remark* 7. An interesting consequence of Theorem 3 is that the optimality of (2) can be implied by the principle of least-squares approximations based on Corollary 1. This is a crucial result since the derivation and structure of (2) does not suggest its optimality with $U \succ 0$ is related to the least-squares approximations.

A particular case of (18) with fewer matrix variables is presented in the following corollary, where the largest lower bound of the inequality is also identical to (5).

**Corollary 2.** *Choosing the same $\varpi(\cdot)$, $U \succ 0$ and $\boldsymbol{f}(\cdot)$ as in Theorems 1 and 2 with $\rho \in \mathbb{N}$, we have*

$$\int_{\mathcal{K}} \varpi(\tau)\boldsymbol{x}^\top(\tau)U\boldsymbol{x}(\tau)\mathrm{d}\tau \geq \boldsymbol{z}^\top \left[\mathbf{Sy}\left(\Upsilon^\top \widehat{X}\right) - \widehat{X}^\top \left(\mathsf{F}^{-1} \otimes U^{-1}\right) \widehat{X}\right] \boldsymbol{z} \tag{25}$$

*with $\widehat{X} = \mathbf{Col}_{i=1}^{d} X_i \in \mathbb{R}^{dn \times \rho n}$ and $\mathsf{F}^{-1} = \int_{\mathcal{K}} \varpi(\tau)\boldsymbol{f}(\tau)\boldsymbol{f}^\top(\tau)\mathrm{d}\tau$, where $\Upsilon$ and $\boldsymbol{z}$ satisfy* (17). *Moreover,* (25) *is a special case of* (18) *and the largest lower bound of* (25) *is identical to* (5) *which occurs with $\widehat{X} = (\mathsf{F} \otimes U)\Upsilon$ for* (25).

*Proof.* Let the same $\varpi(\cdot)$, $U \succ 0$ and $\boldsymbol{f}(\cdot)$ in Theorems 1 and 2 to be given. Now let $Y = X^\top U^{-1} X$, we have $\mathsf{W} = \widehat{X}^\top \left(\mathsf{F}^{-1} \otimes U^{-1}\right) \widehat{X}$ considering (22), with which the inequality in (18) becomes (25). Moreover, one can conclude by the application of Lemma 2 that the largest lower bound of (25) is attained with $\widehat{X} = (\mathsf{F} \otimes U)\Upsilon$ where (25) becomes (5). ■

*Remark* 8. Let $\mathcal{K} = [a, b]$, $\varpi(\cdot) = 1$ and $\boldsymbol{f}(\tau)$ to contain the Legendre polynomials over $[a, b]$, then Lemma 1 in [27] can be obtained from the corresponding (25) with appropriate $\boldsymbol{z}$ and $\Upsilon$ using the substitution $\boldsymbol{x}(\cdot) \leftarrow \dot{\boldsymbol{x}}(\cdot)$. Now consider the fact that the left hand of the inequality (9) in [27] can be rewritten as a one fold integral with a weight function by using the Cauchy formula for repeated integrations. Let $\mathcal{K} = [a, b]$ and $\boldsymbol{f}(\tau)$ to contain Jacobi polynomials associated with $\varpi(\tau) = (\tau - a)^m$ over $[a, b]$, then [27, eq.(9)] can be obtained by the corresponding (25) with appropriate $\boldsymbol{z}$ and $\Upsilon$ using the substitution $\boldsymbol{x}(\cdot) \leftarrow \dot{\boldsymbol{x}}(\cdot)$. Finally, since the largest inequality's lower bound of (25) is identical to the lower bound in (5) for the same $\varpi(\cdot)$, $\boldsymbol{f}(\cdot)$ and $U \succ 0$, it also indicates that equivalence relations in terms of inequality's lower bounds can be established among the inequalities in [11], [27].

*Remark* 9. Let $\mathcal{K} = [a, b]$, $\varpi(\cdot) = 1$ and $\boldsymbol{f}(\tau)$ to contain the Legendre polynomials over $[a, b]$, then the conclusion of Theorem 1 in [14] can be obtained from Theorem 3 with appropriate $\Upsilon$ and $\boldsymbol{z}$ considering the substitution $\boldsymbol{x}(\cdot) \leftarrow \dot{\boldsymbol{x}}(\cdot)$. As we have proved that the lower bound in (5) corresponds to the largest inequality's lower bound of (25), thus the largest inequality's lower bounds of (18) and (25) are identical. Consequently, it is possible to show that equivalent relations[1] in terms of inequality's lower bounds can be established between the inequalities in [11], [14], [27] given what we have presented in Remark 8.

*Remark* 10. Since the inequality's lower bound of (5) with $U \succ 0$ corresponds to the largest lower bound in (25), thus the optimality of (25) is inferred by the optimality of (5) with $U \succ 0$, which is guaranteed by the principle of least-squares approximations as we have proved in Corollary 1.

The conclusion in Theorem 3 is very important in terms of understanding the relationship among (5), (18) and (25). Since all these three inequalities are essentially equivalent in terms of their lower bounds with the same $\varpi(\cdot)$, $\boldsymbol{f}(\cdot)$ and $U \succ 0$, hence if a particular case of one of the three inequalities is derived then the corresponding two 'equivalent' inequalities can be immediately constructed.

*Remark* 11. Because of the result in Theorem 3 together with what we had explained in Remark 5, one can conclude that an invertible linear transformation $G$ acting on $\boldsymbol{f}(\cdot)$, namely, $G\boldsymbol{f}(\cdot) \to \boldsymbol{f}(\cdot)$, gives an equivalent inequality's lower bound as the one in (18).

## IV. STABILITY ANALYSIS OF A COUPLED DIFFERENTIAL DIFFERENCE SYSTEM WITH A DISTRIBUTED DELAY

The proposed inequalities could be applied to various contexts such as the stability analysis and control of infinite-dimensional systems such as delay [1], [2] and ODE-PDE coupled system [5], [6]. To show an application in this paper, we derive two stability conditions in this section for a linear CDDS with a distributed delay via the application of (5) and (25). We also show that the resulting stability conditions are equivalent where their solvability remains unchanged with respect to a matrix parameter in the Liapunov-Krasovskii functional.

Consider a linear CDDS of the form

$$
\begin{aligned}
&\dot{\boldsymbol{x}}(t) = A_1\boldsymbol{x}(t) + A_2\boldsymbol{y}(t - r) + \int_{-r}^{0} \widetilde{A}_3(\tau)\boldsymbol{y}(t + \tau)\mathsf{d}\tau \\
&\boldsymbol{y}(t) = A_4\boldsymbol{x}(t) + A_5\boldsymbol{y}(t - r), \quad t \geq t_0 \\
&\boldsymbol{x}(t_0) = \boldsymbol{\xi}, \quad \forall\theta \in [-r, 0], \quad \boldsymbol{y}(t_0 + \theta) = \boldsymbol{\phi}(\theta)
\end{aligned}
\tag{26}
$$

[1]The equivalence relations here are understood by considering the structure of inequalities irrespective of using $\boldsymbol{x}(\cdot)$ or $\dot{\boldsymbol{x}}(\cdot)$.

where $t_0 \in \mathbb{R}$ and $\boldsymbol{\xi} \in \mathbb{R}^n$ and $\boldsymbol{\phi}(\cdot) \in \widehat{\mathbb{C}}([-r,0) ⨟ \mathbb{R}^\nu)$. The notation $\widehat{\mathbb{C}}([-r,0) ⨟ \mathbb{R}^n)$ stands for the space of bounded right piecewise continuous functions endowed with $\|\boldsymbol{\phi}(\cdot)\|_\infty = \sup_{\tau \in \mathcal{X}} \|\boldsymbol{\phi}(\tau)\|_2$. We also assume that $\rho(A_5) < 1$ which ensures the input to state stability (See [17]) of $\boldsymbol{y}(t) = A_4\boldsymbol{x}(t) + A_5\boldsymbol{y}(t-r)$, where $\rho(A_5)$ is the spectral radius of $A_5$. Since $\rho(A_5) < 1$ is independent from $r$, thus this condition ensures the input to state stability of $\boldsymbol{y}(t) = A_4\boldsymbol{x}(t) + A_5\boldsymbol{y}(t-r)$ for all $r > 0$. Finally, $\widetilde{A}_3(\tau)$ in (26) satisfies the following assumption.

**Assumption 1.** There exist $\mathbf{Col}_{i=1}^d f_i(\tau) = \boldsymbol{f}(\cdot) \in \mathbb{C}^1(\mathbb{R} ⨟ \mathbb{R}^d)$ and $A_3 \in \mathbb{R}^{n\times\nu d}$ with $d \in \mathbb{N}$, which satisfy

$$\int_{-r}^{0} \boldsymbol{f}(\tau)\boldsymbol{f}^\top(\tau)\mathsf{d}\tau \succ 0 \tag{27}$$

$$\exists M \in \mathbb{R}^{d\times d}, \quad \frac{\mathsf{d}\boldsymbol{f}(\tau)}{\mathsf{d}\tau} = M\boldsymbol{f}(\tau), \tag{28}$$

and $\forall \tau \in [-r,0], \widetilde{A}_3(\tau) = A_3F(\tau) \in \mathbb{R}^{n\times\nu}$, where $F(\tau) = \boldsymbol{f}(\tau) \otimes I_\nu \in \mathbb{R}^{\nu d\times\nu}$.

*Remark* 12. The functions in $\boldsymbol{f}(\cdot)$ in Assumption 1 are the solutions of homogeneous differential equations. Many models of delay systems are encompassed by (26), which is the main reason why (26) is chosen as the model to be analyzed in this section. Specifically, see the examples in [2], [17] and the references therein.

The following Liapunov-Krasovskii stability criteria is employed in this section to construct stability conditions for (26).

**Lemma 4.** *Given $r > 0$, the system in* (26) *is globally uniformly asymptotically stable at its origin, if there exist $\epsilon_1; \epsilon_2; \epsilon_3 > 0$ and a differentiable functional $v : \mathbb{R}^n \times \widehat{\mathbb{C}}([-r,0) ⨟ \mathbb{R}^\nu) \to \mathbb{T}$ such that $v(\mathbf{0}_n, \mathbf{0}_\nu) = 0$ and*

$$\epsilon_1\|\boldsymbol{\xi}\|_2^2 \le v(\boldsymbol{\xi}, \boldsymbol{\phi}(\cdot)) \le \epsilon_2 \left(\|\boldsymbol{\xi}\|_2 \vee \|\boldsymbol{\phi}(\cdot)\|_\infty\right)^2 \tag{29}$$

$$\left.\frac{\mathsf{d}^+}{\mathsf{d}t} v(\boldsymbol{x}(t), \mathbf{y}_t(\cdot))\right|_{t=t_0, \boldsymbol{x}(t_0)=\boldsymbol{\xi}, \mathbf{y}_{t_0}(\cdot)=\boldsymbol{\phi}(\cdot)} \le -\epsilon_3\|\boldsymbol{\xi}\|_2^2 \tag{30}$$

*for any $\boldsymbol{\xi} \in \mathbb{R}^n$ and $\boldsymbol{\phi}(\cdot) \in \widehat{\mathbb{C}}([-r,0) ⨟ \mathbb{R}^\nu)$ in* (26)*, where $t_0 \in \mathbb{R}$ and $\frac{\mathsf{d}^+}{\mathsf{d}x} f(x) = \operatorname{limsup}_{\eta\downarrow 0} \frac{f(x+\eta)-f(x)}{\eta}$. Furthermore, $\mathbf{y}_t(\cdot)$ in* (30) *is defined by the equality $\forall t \ge t_0$, $\forall \theta \in [-r,0)$, $\mathbf{y}_t(\theta) = \boldsymbol{y}(t+\theta)$ where $\boldsymbol{x}(t)$ and $\boldsymbol{y}(t)$ satisfying* (26).

*Proof.* Let $u(\cdot), v(\cdot), w(\cdot)$ in Theorem 3 of [17] be quadratic functions with the multiplier factors $\epsilon_1; \epsilon_2; \epsilon_3 > 0$. Since (26) is a particular case of the general system considered in Theorem 3 of [17], then Lemma 4 is obtained. ∎

To analyze the stability of the origin of (26), consider the following parameterized Krasovskii functional

$$v(\boldsymbol{\xi}, \boldsymbol{\phi}(\cdot)) := \begin{bmatrix} \boldsymbol{\xi} \\ \int_{-r}^0 \widehat{G}(\tau)\boldsymbol{\phi}(\tau)\mathsf{d}\tau \end{bmatrix}^\top \widehat{P} \begin{bmatrix} \boldsymbol{\xi} \\ \int_{-r}^0 \widehat{G}(\tau)\boldsymbol{\phi}(\tau)\mathsf{d}\tau \end{bmatrix} + \int_{-r}^0 \boldsymbol{\phi}^\top(\tau) \left[S + (\tau + r)U\right] \boldsymbol{\phi}(\tau)\mathsf{d}\tau \tag{31}$$

with $\widehat{G}(\tau) = \boldsymbol{g}(\tau) \otimes I_\nu$, $\boldsymbol{g}(\tau) = G\boldsymbol{f}(\tau)$, $G \in \mathbb{R}_{[d]}^{d\times d}$ and $\boldsymbol{f}(\cdot)$ in Assumption 1, where $\boldsymbol{\xi} \in \mathbb{R}^n$, $\boldsymbol{\phi}(\cdot) \in \widehat{\mathbb{C}}([-r,0) ⨟ \mathbb{R}^\nu)$ in (31) are the initial conditions in (26), and $\widehat{P} \in \mathbb{S}^{n+d\nu}$ and $S; U \in \mathbb{S}^\nu$ are unknown parameters to be determined. Note that $\widehat{G}(\tau)$ can be rewritten as $\widehat{G}(\tau) = \boldsymbol{g}(\tau) \otimes I_\nu = G\boldsymbol{f}(\tau) \otimes I_\nu = (G \otimes I_\nu)F(\tau)$ with $F(\tau) := \boldsymbol{f}(\tau) \otimes I_\nu$ based on the property of the Kronecker product in (1). Note that also (31) can be regarded as a parameterized version of the complete Liapunov-Krasovskii functional proposed in [17].

We will show in the following theorem that the solvability of the resulting stability conditions remain unchanged for any value of $G \in \mathbb{R}_{[d]}^{d\times d}$ in (31) when (5) or (25) is applied for the derivation.

**Theorem 4.** *Given $G \in \mathbb{R}^{d\times d}_{[d]}$ and $\boldsymbol{f}(\cdot)$, $M$ in* (28) *with $\mathsf{F}^{-1} = \int_{-r}^{0} \boldsymbol{f}(\tau)\boldsymbol{f}^\top(\tau)\mathsf{d}\tau$, then the origin of* (26) *under Assumption 1 is globally uniformly asymptotically stable if there exists $\widehat{P} \in \mathbb{S}^{n+d\nu}$ and $S;U \in \mathbb{S}^{\nu}$ such that*

$$\widehat{P} + \left[\mathsf{O}_n \oplus \left([*]\mathsf{F}G^{-1} \otimes S\right)\right] \succ 0 \tag{32}$$

$$S \succ 0, \;\; U \succ 0, \;\; \mathbf{\Phi}_1 \prec 0 \tag{33}$$

*hold, or equivalently if there exist $\widehat{P} \in \mathbb{S}^{n+d\nu}$, $S;U \in \mathbb{S}^{\nu}$ and $X_1; X_2 \in \mathbb{S}^{d\nu}$ such that*

$$\begin{bmatrix} \widehat{P} + [\mathsf{O}_n \oplus 2X_1] & \begin{bmatrix}\mathsf{O}_{d\nu\times n} & X_1\end{bmatrix}^\top \\ * & [*]\mathsf{F}G^{-1} \otimes S \end{bmatrix} \succ 0 \tag{34}$$

$$S \succ 0, \;\; U \succ 0, \;\; \begin{bmatrix} \mathbf{\Phi}_2 & \begin{bmatrix}\mathsf{O}_{d\nu\times(\nu+n)} & X_2\end{bmatrix}^\top \\ * & -[*]\mathsf{F}G^{-1} \otimes U \end{bmatrix} \prec 0 \tag{35}$$

*hold, where*

$$\begin{aligned} \mathbf{\Phi}_1 = \mathbf{Sy}\left(H\widehat{P}\begin{bmatrix}\mathbf{A}^\top & \mathbf{G}^\top\end{bmatrix}^\top\right) + \Gamma^\top (S + rU)\,\Gamma \\ - \left(\mathsf{O}_n \oplus S \oplus \left([*]\mathsf{F}G^{-1} \otimes U\right)\right) \end{aligned} \tag{36}$$

$$\begin{aligned} \mathbf{\Phi}_2 = \mathbf{Sy}\left(H\widehat{P}\begin{bmatrix}\mathbf{A}^\top & \mathbf{G}^\top\end{bmatrix}^\top\right) + \Gamma^\top (S + rU)\,\Gamma \\ - \left(\mathsf{O}_n \oplus S \oplus 2X_2\right) \end{aligned} \tag{37}$$

$$H = \begin{bmatrix} I_n & \mathsf{O}_{n\times d\nu} \\ \mathsf{O}_{\nu\times n} & \mathsf{O}_{\nu\times d\nu} \\ \mathsf{O}_{d\nu\times n} & I_{d\nu} \end{bmatrix}, \;\; \Gamma := \begin{bmatrix} A_4 & A_5 & \mathsf{O}_{\nu\times d\nu}\end{bmatrix} \tag{38}$$

$$\mathbf{A} = \begin{bmatrix} A_1 & A_2 & A_3(G^{-1}\otimes I_\nu)\end{bmatrix} \tag{39}$$

$$\mathbf{G} = \begin{bmatrix} \widehat{G}(0)A_4 & \widehat{G}(0)A_5 - \widehat{G}(-r) & -\widehat{M}\end{bmatrix} \tag{40}$$

*with $\widehat{G}(0) = (G \otimes I_\nu)F(0) = G\boldsymbol{f}(0) \otimes I_\nu$ and $\widehat{G}(-r) = (G \otimes I_\nu)F(-r) = G\boldsymbol{f}(-r) \otimes I_\nu$ and $\widehat{M} = (G \otimes I_\nu)(M \otimes I_\nu)(G^{-1} \otimes I_\nu) = GMG^{-1} \otimes I_\nu$. Finally, the value of $G \in \mathbb{R}^{d\times d}_{[d]}$ does not change the solvability of the matrix inequalities in* (32)–(35).

*Proof.* Let $G \in \mathbb{R}^{d\times d}_{[d]}$ and $\boldsymbol{f}(\cdot)$ with $M$ in Assumption 1 be given. Given the fact that the eigenvalues of $S+(\tau+r)U$, $\tau \in [-r, 0]$ are bounded and $\widehat{G}(\tau) = (G \otimes I_n)F(\tau)$, it is obvious to see that (31) satisfies the following property that there exist $\lambda; \eta > 0$ such that

$$\begin{aligned} v(\boldsymbol{\xi}, \boldsymbol{\phi}(\cdot)) &\le \begin{bmatrix} \boldsymbol{\xi} \\ \int_{-r}^{0} F(\tau)\boldsymbol{\phi}(\tau)\mathsf{d}\tau \end{bmatrix}^\top \lambda \begin{bmatrix} \boldsymbol{\xi} \\ \int_{-r}^{0} F(\tau)\boldsymbol{\phi}(\tau)\mathsf{d}\tau \end{bmatrix} + \int_{-r}^{0} \boldsymbol{\phi}^\top(\tau)\lambda\boldsymbol{\phi}(\tau)\mathsf{d}\tau \\ &\le \lambda\|\boldsymbol{\xi}\|_2^2 + \int_{-r}^{0} \boldsymbol{\phi}^\top(\tau)F^\top(\tau)\mathsf{d}\tau\lambda \int_{-r}^{0} F(\tau)\boldsymbol{\phi}(\tau)\mathsf{d}\tau + \lambda r\|\boldsymbol{\phi}(\cdot)\|_\infty^2 \\ &\le \lambda\|\boldsymbol{\xi}\|_2^2 + \lambda r\|\boldsymbol{\phi}(\cdot)\|_\infty^2 + \int_{-r}^{0} \boldsymbol{\phi}^\top(\tau)F^\top(\tau)\mathsf{d}\tau\,(\eta\mathsf{F} \otimes I_n) \int_{-r}^{0} F(\tau)\boldsymbol{\phi}(\tau)\mathsf{d}\tau \\ &\le \lambda\|\boldsymbol{\xi}\|_2^2 + \lambda r\|\boldsymbol{\phi}(\cdot)\|_\infty^2 + \int_{-r}^{0} \boldsymbol{\phi}^\top(\tau)\eta\boldsymbol{\phi}(\tau)\mathsf{d}\tau \\ &\le \lambda\|\boldsymbol{\xi}\|_2^2 + (\lambda r + \eta r)\,\|\boldsymbol{\phi}(\cdot)\|_\infty^2 \le (\lambda + \lambda r + \eta r)\left(\|\boldsymbol{\xi}\|_2^2 + \|\boldsymbol{\phi}(\cdot)\|_\infty^2\right) \end{aligned}$$

$$\leq 2\left(\lambda+\lambda r+\eta r\right)\left[\max\left(\|\boldsymbol{\xi}\|_2, \|\boldsymbol{\phi}(\cdot)\|_\infty\right)\right]^2 \quad (41)$$

for any $\boldsymbol{\xi} \in \mathbb{R}^n$ and $\boldsymbol{\phi}(\cdot) \in \widehat{\mathbb{C}}\left([-r_2, 0) ⨾ \mathbb{R}^\nu\right)$ in (26), where (41) is derived via the property of quadratic forms: $\forall X \in \mathbb{S}^n$, $\exists \lambda > 0 : \forall \mathbf{x} \in \mathbb{R}^n \setminus \{\mathbf{0}\}$, $\mathbf{x}^\top(\lambda I_n - X)\,\mathbf{x} > 0$ together with the application of (5) with $\boldsymbol{f}(\cdot)$ in (1). Consequently, (41) shows that (31) satisfies the upper bound property in (29).

Now apply (5) with $\varpi(\tau) = 1$ to $\int_{-r}^{0} \boldsymbol{\phi}^\top(\tau)S\boldsymbol{\phi}(\tau)\mathsf{d}\tau$ in (31) given $S \succ 0$ and $\boldsymbol{f}(\cdot)$ in Assumption 1 and the fact that $\boldsymbol{\phi}(\cdot) \in \widehat{\mathbb{C}}([-r_2, 0) ⨾ \mathbb{R}^\nu) \subset \mathbb{L}^2([-r_2, 0) ⨾ \mathbb{R}^\nu)$. It yields

$$\int_{-r}^{0} \boldsymbol{\phi}^\top(\tau)S\boldsymbol{\phi}(\tau)\mathsf{d}\tau \geq \int_{-r}^{0} \boldsymbol{\phi}^\top(\tau)\widehat{G}^\top(\tau)\mathsf{d}\tau \left[[*]\mathsf{F}G^{-1} \otimes S\right] \int_{-r}^{0} \widehat{G}(\tau)\boldsymbol{\phi}(\tau)\mathsf{d}\tau \quad (42)$$

for any $\boldsymbol{\xi} \in \mathbb{R}^n$ and $\boldsymbol{\phi}(\cdot) \in \widehat{\mathbb{C}}\left([-r_2, 0) ⨾ \mathbb{R}^\nu\right)$ in (26).

Moreover, by applying (42) to the functional in (31), we can conclude that if (32) is feasible, then it infers the existence of $\epsilon_1; \epsilon_2 > 0$ and the functional in (31) satisfying (29) given what we have shown in (41). On the other hand, considering the property of congruence transformations with the fact that $G \in \mathbb{R}_{[d]}^{d\times d}$, one can conclude that (32) holds if and only if

$$\begin{aligned}&\left[I_n \oplus \left(G^\top \otimes I_\nu\right)\right] \widehat{P} \left[I_n \oplus (G \otimes I_\nu)\right] + \\ &\qquad \left[I_n \oplus \left(G^\top \otimes I_\nu\right)\right] \left[\mathsf{O}_n \oplus \left([*]\mathsf{F}G^{-1} \otimes S\right)\right] \left[I_n \oplus (G \otimes I_\nu)\right] \\ &\qquad\qquad = P + \left[\mathsf{O}_n \oplus (\mathsf{F} \otimes S)\right] \succ 0 \end{aligned} \quad (43)$$

with $P = \left[I_n \oplus \left(G^\top \otimes I_\nu\right)\right] \widehat{P} \left[I_n \oplus (G \otimes I_\nu)\right]$. By treating $P$ as a new variable, it shows that the solvability of the last matrix inequality in (43), namely (32), is not affected by the values of $G \in \mathbb{R}_{[d]}^{d\times d}$.

Now we start to prove that if (32) or (35) is feasible then there exist $\epsilon_3 > 0$ and (31) satisfying (30). Differentiate $v(\boldsymbol{x}(t), \mathbf{y}_t(\cdot))$ along the trajectory of (26) at $t = t_0$ and consider the relation

$$\begin{aligned}\left.\frac{\mathsf{d}}{\mathsf{d}t}\int_{-r}^{0} \widehat{G}(\tau)\boldsymbol{y}(t+\tau)\mathsf{d}\tau\right|_{t=t_0} &= (G \otimes I_\nu)F(0)\boldsymbol{\phi}(0) - (G \otimes I_\nu)F(-r)\boldsymbol{\phi}(-r) - \widehat{M}\int_{-r}^{0} \left(G\boldsymbol{f}(\tau) \otimes I_\nu\right)\boldsymbol{\phi}(\tau)\mathsf{d}\tau \\ &= \widehat{G}(0)A_4\boldsymbol{\xi} + \left[\widehat{G}(0)A_5 - \widehat{G}(-r)\right]\boldsymbol{\phi}(-r) - \widehat{M}\int_{-r}^{0} \widehat{G}(\tau)\boldsymbol{\phi}(\tau)\mathsf{d}\tau, \end{aligned} \quad (44)$$

where $\widehat{M} = (G \otimes I_\nu)(M \otimes I_\nu)(G^{-1} \otimes I_\nu)$ and (44) is obtained via (28). Then we have

$$\begin{aligned}\left.\frac{\mathsf{d}^+}{\mathsf{d}t}v(\boldsymbol{x}(t), \mathbf{y}_t(\cdot))\right|_{t=t_0, \boldsymbol{x}(t_0)=\boldsymbol{\xi}, \mathbf{y}_{t_0}(\cdot)=\boldsymbol{\phi}(\cdot)} &= \boldsymbol{\chi}^\top \mathbf{Sy}\left(H\widehat{P}\begin{bmatrix}\mathbf{A}\\ \mathbf{G}\end{bmatrix}\right)\boldsymbol{\chi} - \int_{-r}^{0} \boldsymbol{\phi}^\top(\tau)U\boldsymbol{\phi}(\tau)\mathsf{d}\tau \\ &\quad + \boldsymbol{\chi}^\top\left[\Gamma^\top(S + rU)\,\Gamma - (\mathsf{O}_n \oplus S \oplus \mathsf{O}_{d\nu})\right]\boldsymbol{\chi} \end{aligned} \quad (45)$$

where $H$, $\mathbf{A}$, $\mathbf{G}$ and $\Gamma$ have been given in (38)–(40) and

$$\boldsymbol{\chi} := \mathbf{Col}\left(\boldsymbol{\xi},\ \boldsymbol{\phi}(-r),\ \int_{-r}^{0} \widehat{G}(\tau)\boldsymbol{\phi}(\tau)\mathsf{d}\tau\right). \quad (46)$$

Given $U \succ 0$ in (32), apply the inequality (5) with $\varpi(\tau) = 1$ to $\int_{-r}^{0} \boldsymbol{\phi}^\top(\tau)U\boldsymbol{\phi}(\tau)\mathsf{d}\tau$ in (45). It yields

$$\int_{-r}^{0} \boldsymbol{\phi}^\top(\tau)U\boldsymbol{\phi}(\tau)\mathsf{d}\tau \geq \int_{-r}^{0} \boldsymbol{\phi}^\top(\tau)\widehat{G}^\top(\tau)\mathsf{d}\tau \left[[*]\mathsf{F}G^{-1} \otimes U\right] \int_{-r}^{0} \widehat{G}(\tau)\boldsymbol{\phi}(\tau)\mathsf{d}\tau \quad (47)$$

for any $\boldsymbol{\xi} \in \mathbb{R}^n$ and $\boldsymbol{\phi}(\cdot) \in \widehat{\mathbb{C}}\left([-r_2, 0) ⨾ \mathbb{R}^\nu\right)$ in (26). By using (47) to (45) with $U \succ 0$, we have

$$\left.\frac{\mathsf{d}^+}{\mathsf{d}t}v(\boldsymbol{x}(t), \mathbf{y}_t(\cdot))\right|_{t=t_0, \boldsymbol{x}(t_0)=\boldsymbol{\xi}, \mathbf{y}_{t_0}(\cdot)=\boldsymbol{\phi}(\cdot)} \leq \boldsymbol{\chi}^\top\boldsymbol{\Phi}_1\boldsymbol{\chi} \quad (48)$$

where $\mathbf{\Phi}_1$ is given in (36). By (48) and (46), it is easy to see that the feasible solutions of (33) infer the existence of $\epsilon_3 > 0$ and (31) satisfying (30).

Now considering the property of congruence transformations with the fact that $G \in \mathbb{R}_{[d]}^{d\times d}$, it is true that

$$\mathbf{\Phi}_1 \prec 0 \iff [*]\,\mathbf{\Phi}_1\left[I_{n+\nu} \oplus (G \otimes I_\nu)\right] = \mathbf{\Theta} := \mathbf{Sy}\left(HP\Psi\right) + \Gamma^\top(S + rU)\,\Gamma - \left[\mathsf{O}_n \oplus S \oplus (\mathsf{F} \otimes U)\right] \prec 0 \quad (49)$$

with $P = \left[I_n \oplus \left(G^\top \otimes I_\nu\right)\right] \widehat{P}\left[I_n \oplus (G \otimes I_\nu)\right]$ and $H$ in (38) and

$$\Psi = \begin{bmatrix} A_1 & A_2 & A_3 \\ F(0)A_4 & F(0)A_5 - F(-r) & -M \otimes I_\nu \end{bmatrix} \quad (50)$$

which can be derived via (1). By (49) and (50), it is clear to see that the solvability of $\mathbf{\Phi}_1 \prec 0$ in (49) is not affected the values of $G \in \mathbb{R}_{[n]}^{n\times n}$ by treating $P$ as a new variable similar to the situation in (43).

Finally, let $\Upsilon = I_{d\nu}$ and $\widehat{X} = X_1 \in \mathbb{S}^{d\nu}$; $\widehat{X} = X_2 \in \mathbb{S}^{d\nu}$ for (25), respectively, then one can apply the inequality in (25) with $\boldsymbol{f}(\cdot)$ in Assumption 1 for the steps at (42) and (47) to construct the constraints in (34)–(35) via the Schur complement. By Theorem 3, one can conclude that (34)–(35) is equivalent to (32)–(33). Since the feasibility of (32)–(33) is not affected by the value of $G \in \mathbb{R}_{[n]}^{n\times n}$, thus the feasibility of (32)–(33) remains unchanged with respect to the value of $G$. ■

*Remark* 13. Although the values of $G \in \mathbb{R}_{[n]}^{n\times n}$ do not affect the solvability of (32)–(35), it can be still beneficial to use orthonormal functions $G\boldsymbol{f}(\tau)$ with $G^{-2} = \int_{\mathcal{K}} \varpi(\tau)\boldsymbol{f}(\tau)\boldsymbol{f}^\top(\tau)\mathsf{d}\tau$ in Theorem 4, as it gives $[*]\mathsf{F}G^{-1} = I_d$ which makes the relevant diagonal blocks in (32)–(35) more regular for numerical calculations.

## V. Numerical example

The numerical example in this section is tested in Matlab with Yalmip [28] for the optimization modeling and Mosek 8 [29] as the numerical solver for semidefinite programming problems. Consider the following linear CDDS

$$\begin{aligned} \dot{x}(t) &= 0.35x(t) + 0.035y(t-r) - \int_{-r}^{0} 5\cos(12\tau)y(t+\tau)\mathsf{d}\tau \\ y(t) &= x(t) + 0.1y(t-r) \end{aligned} \quad (51)$$

which corresponds to $A_1 = 0.35$, $A_2 = 0.035$, $\widetilde{A}_3(\tau) = -5\cos(12\tau)$, $A_4 = 1$ and $A_5 = 0.1$ with $n = 1$ in (26). Now let $\boldsymbol{f}(\cdot)$ and $M$ in (28) be

$$\boldsymbol{f}(\tau) = \begin{bmatrix} 1 \\ \sin(12\tau) \\ \cos(12\tau) \end{bmatrix}, \quad M = \begin{bmatrix} 0 & 0 & 0 \\ 0 & 0 & 12 \\ 0 & -12 & 0 \end{bmatrix} \quad (52)$$

with $A_3 = \begin{bmatrix} 0 & 0 & -5 \end{bmatrix}$ which satisfies Assumption 1.

We will calculate the delay margins of (51) via (32)–(33) and (34)–(35) with (52) and different values of $G$ to show that the numerical results are not affected by the different values of an invertible $G$.

Now apply (32)–(33) with $G = I_3$ and (52) to the system in (51), we obtain the delay boundaries of the detectable stable points of (52) as

$$[0.084, 0.178], \quad [0.607, 0.702], \quad [1.131, 1.225], \quad [1.654, 1.749], \quad [2.178, 2.273] \quad (53)$$

by using a testing vector $r = (1:2500)/1000$, where at each point of $r$ the corresponding optimization program has **12** decision variables. By using (32)–(33) with $G = \begin{bmatrix} 1 & 0.5 & 0.2 \\ 0 & 2 & -1 \\ 0 & 0 & 2 \end{bmatrix}$ and (52) to the system in (51), the **same** results in (53) can be obtained where for each point of $r$ the **same** number of decision variables is required as for $G = I_3$. On the other hand, the **same** boundaries of (53) can be obtained by using the stability condition (34)–(35) with (52) and $G = I_3$ or $G = \begin{bmatrix} 1 & 0.5 & 0.2 \\ 0 & 2 & -1 \\ 0 & 0 & 2 \end{bmatrix}$ to the system in (51), which requires **18** decision variables for each point of $r = (1:2500)/1000$. The above results concerning the stability of (52) are consistent with what have been proved in Theorem 4 that (32)–(33) and (34)–(35) are equivalent and their solvability is not affected the values of $G \in \mathbb{R}^{d\times d}_{[d]}$.

Now we want to show the impact of considering different $\boldsymbol{f}(\cdot)$. Here we modify the conditions (32)–(33) with $G = I_d$ and $\boldsymbol{f}(\tau) = \boldsymbol{\ell}_d(\tau)$ in accordance to the polynomials approximation approach in [30], where $\boldsymbol{\ell}_d(\tau) = \mathbf{Col}_{i=0}^{d}\, \ell_i(\tau)$ contains Legendre polynomials $\{\ell_i(\tau)\}_{i=0}^{d}$ over $[-r, 0]$ up to degree $d \in \mathbb{N}_0$. It occurs that the modified condition requires $d = 22$ with **302** decision variables for each point of $r = (1:2500)/1000$ to detect the boundaries in (53), which is in sharp contrast with the numbers of variables required by (32)–(33) and (34)–(35) with (52). The above results show the advantage of using (52) with our methods over the polynomials approximation approach in [30].

## VI. Conclusion

In this note, two general classes of integral inequalities have been proposed whose integral kernels can be any function in $\mathbb{L}^2_{\varpi}(\mathcal{K}⨾\mathbb{R})$ with a general weight function $\varpi(\tau)$ in (3). Consequently, the proposed inequalities generalize many existing integral inequalities in the existing literature. It has been demonstrated in this note that the optimality of (5) and (18) with $U \succ 0$ is guaranteed by the principle of least-squares approximation based on the results in Corollary 1 and Theorem 3. Moreover, the result concerning (5) and (18) established in Theorem 3 can lead to the conclusion that the lower bounds of many existing quadratic integral inequalities are essentially equivalent. For a specific application, the inequalities presented in this note have been utilized to derive equivalent stability conditions for a linear CDDS with a distributed delay. Finally, the proposed inequalities have great potential to be applied in wider contexts such as the stability analysis of PDE-related systems or sampled-data systems or other types of infinite-dimensional systems whenever the situations are suitable.